\documentclass[a4paper]{article}
\usepackage[latin1]{inputenc}
\usepackage{amsmath,amsthm,amsfonts}
\usepackage{xspace,url}
\usepackage{fullpage}
\usepackage[pdftex,colorlinks=true]{hyperref}
\usepackage[numbers]{natbib}
\newcommand{\qconjugate}[1]{\overline{#1}} 
\newcommand{\cconjugate}[1]{{#1^\star}}    
\renewcommand{\v}[1]{\boldsymbol{#1}} 
\renewcommand{\i}{\ensuremath{\v{i}}\xspace}
\renewcommand{\j}{\ensuremath{\v{j}}\xspace}
\renewcommand{\k}{\ensuremath{\v{k}}\xspace}
\let\oldmu\mu\renewcommand{\mu}{\v\oldmu}
\let\oldnu\nu\renewcommand{\nu}{\v\oldnu}
\let\oldxi\xi\renewcommand{\xi}{\v\oldxi}

\newcommand{\inner}[2]{\ensuremath{\left\langle #1, #2 \right\rangle}\xspace}
\newcommand{\norm}[1]{\ensuremath{\left\|#1\right\|}\xspace}

  \newcommand{\R}{\ensuremath{\mathbb{R}}\xspace}
\renewcommand{\H}{\ensuremath{\mathbb{H}}\xspace}
  \newcommand{\C}{\ensuremath{\mathbb{C}}\xspace}
  \newcommand{\B}{\ensuremath{\mathbb{B}}\xspace}
\newcommand{\I}{\ensuremath{I}\xspace} 
\newtheorem{lemma}{Lemma}
\newtheorem{theorem}{Theorem}
\newtheorem{corollary}{Corollary}
\newtheorem{proposition}{Proposition}
\title{Determination of the biquaternion divisors of zero,
       including the idempotents and nilpotents}
\author{Stephen~J.~Sangwine\thanks{Stephen J. Sangwine is with the
        School of Computer Science and Electronic Engineering,
        University of Essex, Wivenhoe Park, Colchester, CO4 3SQ,
        United Kingdom.
        Email:~\texttt{S.Sangwine@IEEE.org}}
        \and
        Daniel Alfsmann\thanks{Daniel Alfsmann is with the
                               Digital Signal Processing Group (DISPO),
                               Ruhr-Universität Bochum, IC 5/50,
                               44780 Bochum, Germany. 
        Email:~\texttt{alfsmann@nt.ruhr-uni-bochum.de}}\footnotemark[2]
      }
\begin{document}
\maketitle
\begin{abstract}
The biquaternion (complexified quaternion) algebra contains
idempotents (elements whose square remains unchanged) and
nilpotents (elements whose square vanishes).
It also contains divisors of zero (elements with vanishing norm).
The idempotents and nilpotents are subsets of the divisors of zero.
These facts have been reported in the literature, but remain obscure
through not being gathered together using modern notation and terminology.
Explicit formulae for finding all the idempotents, nilpotents and
divisors of zero appear not to be available in the literature, and
we rectify this with the present paper.
Using several different representations for biquaternions, we present
simple formulae for the idempotents, nilpotents and divisors of zero,
and we show that the complex components of a biquaternion
divisor of zero must have a sum of squares that vanishes,
and that this condition is equivalent to two conditions on the inner
product of the real and imaginary parts of the biquaternion, and the
equality of the norms of the real and imaginary parts.
We give numerical examples of nilpotents, idempotents and other divisors of
zero. Finally, we conclude with a statement about the composition
of the set of biquaternion divisors of zero, and its subsets,
the idempotents and the nilpotents.
\end{abstract}
\section{Introduction}
\label{intro}
The complexified quaternions, or biquaternions\footnote{Biquaternions
  was a word coined by Hamilton \cite{Hamiltonpapers:V3:35} and
  \cite[§\,669, p\,664]{Hamilton:1853}. (The word was used 18 years later
  by Clifford \cite{10.1112/plms/s1-4.1.381} for an unrelated type of
  hypercomplex number, which is unfortunate.)},
were discovered by Hamilton himself and published within ten years of the
discovery of quaternions \cite{Hamiltonpapers:V3:35}.
Hamilton discovered that a subset of the biquaternions were divisors of zero
\cite[Lecture VII, §\,672, p\,669]{Hamilton:1853}\footnote{Hamilton did not,
  of course, use the modern terminology \emph{divisor of zero}, but referred
  to the property of a vanishing \emph{tensor}, or what is now called a
  \emph{semi-norm}. (Hamilton's term \emph{tensor} is unrelated to the modern term).}.
He also realised that some divisors of zero had the property of squaring to
yield zero, and thus are solutions of $q^2=0$, which Peirce called \emph{nilpotent} elements
\cite[§\,25, p\,8]{Peirce1882}\cite{Waerden1985}.
As far as we can establish, Hamilton did not realise that some biquaternion
divisors of zero had the property of \emph{idempotency} \cite{Peirce1882,Waerden1985}, that is, squaring
to yield themselves, and thus being solutions of $q(q-1)=0$.
The existence
of biquaternion idempotents and at least one example, has been noted, for example
in \cite[§\,7]{arXiv:physics/0508036} (quoting a paper by Lanczos dating from the 1920s).
Of course, algebraic theorists have established many general properties of
algebras, and the existence of idempotents and nilpotents is widely known
in many algebras, including the biquaternions. However, explicit solutions
for the biquaternion idempotents and nilpotents appear not to have been
worked out thoroughly, and that is the purpose of this paper.

In this paper, we review Hamilton's results on divisors of zero, using modern
terminology, and we solve for the sets of biquaternion idempotents and nilpotents.
The methodology followed is similar to that used in
\cite{arXiv:math.RA/0506190,10.1007/s00006-006-0005-8}
to find the biquaternion square roots of $-1$, but we utilise some different
representations of the biquaternions in our solutions, demonstrating the utility
of these different representations, depending on the algebraic task in hand.

In §\,\ref{notations} we present the notations used in the paper,
and the representations of biquaternions used in the derivations.
In §\,\ref{divisors} we review the biquaternion divisors of zero,
restating Hamilton's results in more modern terms.
Then, in §\,\ref{idempotents} we derive the idempotents
and show that they are a subset of the divisors of zero.
In §\,\ref{nilpotents} we do the same for the nilpotents.
We present examples of divisors of zero, idempotents
and nilpotents in sections \ref{divisors}, \ref{idempotents} and \ref{nilpotents} respectively.

\section{Notations and representations}
\label{notations}
A biquaternion can be represented in many ways, not all of which are discussed here.
We denote the set of reals by \R, the set of complex numbers by \C, the set of (real)
quaternions by \H, and the set of biquaternions (quaternions over the complex field) by \B.

The simplest representation of a biquaternion is the Cartesian representation with
four complex components:
$q = W + X\i + Y\j + Z\k$ where $W, X, Y, Z\in\C$ and \i, \j and \k are the usual
quaternion roots of $-1$.
We need an unambiguous representation for the square root of $-1$ in \C, which is
distinct from the quaternion \i.
In this paper we use the notation \I for this.
\I commutes with the quaternion roots of $-1$,
and hence with all quaternions and biquaternions.

Two alternative representations are used in this paper.
The first is $q = q_r + \I q_i$, in which
$q_r, q_i\in\H$. In this form, the biquaternion is represented as an ordered pair
of quaternions, or a complex number with quaternion real and imaginary parts.
Expressing $q_r = w_r + x_r\i + y_r\j + z_r\k$,
and $q_i = w_i + x_i\i + y_i\j + z_i\k$,
$w_a, x_a, y_a, z_a\in\R, a\in\{r,i\}$
we have $w_r = \Re(W)$, $w_i = \Im(W)$, or $W = w_r + \I w_i$,
and so on.

The Cartesian components $w$ and $W$ are known as \emph{scalar} parts and
the remaining components of a quaternion or biquaternion are known as the
\emph{vector} part.
A quaternion with $w=0$ or a biquaternion with $W=0$ (that is, a zero scalar
part) is called \emph{pure} (terminology due to Hamilton).

In the complex representation $q = q_r + \I q_i$,
we make use of the inner product of the real and imaginary
quaternions, $\inner{q_r}{q_i} = w_r w_i + x_r x_i + y_r y_i + z_r z_i$,
as defined by Porteous \cite[Prop.\,10.8, p\,177]{Porteous:1981}.
The inner product vanishes when the two quaternions are perpendicular in 4-space.

The second representation that we use is $q = A + \xi B$, where $A, B\in\C, \xi\in\B$ and $\xi^2=-1$.
In this form, the biquaternion is represented as a complex number with complex
real and imaginary parts ($A$ and $B$), and $\xi$ is a biquaternion root of $-1$
as defined in \cite{arXiv:math.RA/0506190,10.1007/s00006-006-0005-8}%
\footnote{The conditions required for $\xi$ to be a root of $-1$ are as
follows: $\xi$ must be a \emph{pure} biquaternion, with $W=0$;
writing $\v\alpha=\Re(\xi)$ and $\v\beta=\Im(\xi)$, $\v\alpha\perp\v\beta$ and
$\norm{\v\alpha}-\norm{\v\beta}=1$.}.
In this second representation,
$A$ is equal to $W$ in the first representation -- they are both the
complex scalar part of the biquaternion, and $B$ is the complex modulus of the
vector part.
$\xi$ is referred to as the `axis' of the biquaternion
by generalization from the quaternion case, where the axis defines the direction
in 3-space of the vector part.
In terms of the Cartesian representation $B = \sqrt{X^2 + Y^2 + Z^2}$
and $\xi = (X\i + Y\j + Z\k)/B$.
Of course, we must be careful here that $B\ne0$, which occurs in the case
where $q$ is a nilpotent, as discussed in §\,\ref{nilpotents}. If $q$ is
nilpotent, $A=0$ and since a nilpotent is also a divisor of zero as shown
later in Theorem~\ref{vanishsquare}, if we try to find $B$ we obtain zero.
Therefore, we cannot find $\xi$ by dividing $q$ by $B$.

In what follows we make use of the conjugate of a quaternion,
which we denote by $\qconjugate{q} = w - x\i - y\j - z\k$
and we also use the concept of a \emph{norm}, which for the real quaternions is
simply defined as $\norm{q} = w^2 + x^2 + y^2 + z^2$, where $w, x, y, z\in\R$.
We generalise both definitions to the case where $q\in\B$ as
$\qconjugate{q} = W - X\i - Y\j - Z\k = A - \xi B$
and $\norm{q} = W^2 + X^2 + Y^2 + Z^2$.
The norm may also be computed as $\norm{q} = q\qconjugate{q} =\qconjugate{q}q$ in both
the quaternion and biquaternion cases.
Conventionally a norm is real and positive definite, but in the case of biquaternions,
this convention must be relaxed, because it is possible for the sum of the squares of the four complex
components to vanish, even though individually they are not all zero.
The term \emph{semi-norm} is therefore used instead.
A semi-norm is a generalization of the concept of a norm, with no requirement that the norm be
zero only at the origin \cite[See: \textbf{semi-norm}]{CollinsDictMaths}.

\section{Divisors of zero}
\label{divisors}
It is possible for the semi-norm of a biquaternion \norm{q} to be zero, even though
$q\ne0$. This cannot happen with a real quaternion
(or, in fact, with a purely imaginary
quaternion). The conditions for the semi-norm to vanish were discovered by Hamilton
\cite[Lecture VII, §\,672, p\,669]{Hamilton:1853}. In modern terminology, an element
of an algebra which has a vanishing semi-norm is known as a \emph{divisor of zero}
since it will satisfy the equation $q a = 0$ for some $a$, even though $q\ne 0$ and $a\ne 0$.
We now find the biquaternion divisors of zero in modern notation.
\begin{lemma}\label{normlemma}
Let $q = q_r + \I q_i$ be a non-zero biquaternion with real part $q_r\in\H$
and imaginary part $q_i\in\H$.
The real part of the semi-norm of $q$ is equal to the difference between
the norms of $q_r$ and $q_i$, and
the imaginary part of the semi-norm of $q$ is equal to twice the inner product of $q_r$ and $q_i$.
\end{lemma}
\begin{proof}
We express the semi-norm of $q$ as the product of $q$ with its quaternion conjugate:
$\norm{q} = q\qconjugate{q}$. Writing this explicitly:
\begin{align*}
\norm{q} &= (q_r + \I q_i)\qconjugate{(q_r + \I q_i)} = (q_r + \I q_i)(\qconjugate{q_r} + \I\qconjugate{q_i})\\
\intertext{and multiplying out we get:}
\label{normq}
\norm{q} &= q_r\qconjugate{q_r} - q_i\qconjugate{q_i} + \I(q_r\qconjugate{q_i} + q_i\qconjugate{q_r})
\end{align*}
The real part of \norm{q} (with respect to \I) can be recognised as $\norm{q_r} - \norm{q_i}$
from Coxeter's result \cite[§\,2]{Coxeter:1946}.
The imaginary part of \norm{q} is twice the inner product of $q_r$ and $q_i$ as given by
Porteous \cite[Prop.\,10.8, p\,177]{Porteous:1981}.
\end{proof}
\begin{theorem}
\label{vanish}
Let $q = q_r + \I q_i$ be a non-zero biquaternion with real part $q_r$ and imaginary part $q_i$.
Then the semi-norm of $q$ is zero iff $\norm{q_r} = \norm{q_i}$ and $\inner{q_r}{q_i}=0$,
that is, the real and imaginary parts must have equal norms and a vanishing inner product.
\end{theorem}
\begin{proof}
If $\norm{q} = 0$, the real and imaginary parts must separately be zero. From Lemma \ref{normlemma}
this requires $\norm{q_r}=\norm{q_i}$ and $\inner{q_r}{q_i}=0$.
\end{proof}
Notice that Theorem \ref{vanish} may be satisfied in the case of a pure biquaternion.
We return to this case in Theorem \ref{vanishsquare} below.

\begin{proposition}
\label{normequivalence}
The conditions for the semi-norm to vanish, stated in Theorem~\ref{vanish},
are equivalent to the condition $\norm{q} = W^2 + X^2 + Y^2 + Z^2 = 0$.
\end{proposition}
\begin{proof}
Let $W = w_r + w_i\I$ where $w_r, w_i\in\R$ and similarly for $X$, $Y$ and $Z$.
Then $W^2 = w_r^2 - w_i^2 + 2w_r w_i\I$ and similarly for $X^2$, $Y^2$ and $Z^2$.
Thus:
\[
W^2 + X^2 + Y^2 + Z^2 = 0 \implies
\begin{aligned}
&\! w_r^2 - w_i^2 + 2w_r w_i\I\\
+\,& x_r^2 - x_i^2 + 2x_r x_i\I\\
+\,& y_r^2 - y_i^2 + 2y_r y_i\I\\
+\,& z_r^2 - z_i^2 + 2z_r z_i\I\\
\end{aligned}
= 0
\]
Equating real and imaginary parts, we obtain the simultaneous equations:
\begin{align*}
w_r^2 - w_i^2 + x_r^2 - x_i^2 + y_r^2 - y_i^2 + z_r^2 - z_i^2 &= 0\\
w_r w_i       + x_r x_i       + y_r y_i       + z_r z_i       &= 0
\end{align*}
The first equation may be grouped as $(w_r^2 + x_r^2 + y_r^2 + z_r^2) - (w_i^2 + x_i^2 + y_i^2 + z_i^2) = 0$
which shows that the norms of the real and imaginary parts must be equal, which is the first
condition stated in Theorem~\ref{vanish}.
The second equation shows that the inner product of the real and imaginary parts must be
zero, which is the second condition stated in Theorem~\ref{vanish}.
\end{proof} 
Notice that Proposition~\ref{normequivalence} holds when $W=0$ and the biquaternion is pure.
This is important for the case of nilpotents, discussed in §\,\ref{nilpotents},
Proposition~\ref{nilpotentsaredivisors}.

Examples of biquaternions with vanishing semi-norm may easily be constructed to meet the
constraints given in Theorem \ref{vanish}. For example:
$q = c + \cconjugate{c}\i - c\j - \cconjugate{c}\k$ where $c = 1 + \I$ and the notation
$\cconjugate{c}$ means the complex conjugate, that is $\cconjugate{c} = 1 - \I$.
The semi-norm is:
\begin{align*}
\norm{q} &= c^2 + \cconjugate{c}^2 + (-c)^2  + (-\cconjugate{c})^2\\
         &= (1 + \I)^2 + (1 - \I)^2 + (-1 - \I)^2  + (-1 + \I)^2\\
         &= 1 - 1 + 2\I + 1 - 1 - 2\I + 1 - 1 + 2\I + 1 - 1 - 2\I = 0
\end{align*}
The real part of this biquaternion is $1 + \i - \j - \k$ and the imaginary part is
$1 - \i - \j + \k$, and these clearly have the same norm, namely 4.
That they are perpendicular in 4-space is easily seen by evaluating their scalar
product: $1 × 1 + 1 × (-1) + (-1) × (-1) + (-1) × 1 = 1 - 1 + 1 - 1 = 0$.

A second example is: $p = \sqrt{8} + \I(2\j + 2\k)$,
because the real and imaginary parts are perpendicular
(the scalar product of $\sqrt{8}$ with $2\j+2\k$ is zero), and the real and imaginary
parts have the same norm, $\norm{p_r}=\norm{p_i}=8$.

\section{Idempotents}
\label{idempotents}
An idempotent is a value that squares to give itself. An example of a biquaternion idempotent is
given in \cite{arXiv:physics/0508036} (quoting a paper of Lanczos) but no general case is given.
We present here a derivation of the subset of biquaternions which are idempotent.
\begin{theorem}
\label{idempotent}
Any biquaternion of the form $q = \frac{1}{2}\pm\frac{1}{2}\xi\I$, where $\xi\in\B$ is a root of $-1$
is an idempotent. There are no other idempotents in \B.
\end{theorem}
\begin{proof}
The proof is by construction from an arbitrary biquaternion $q$, represented in the form
$q = A + \xi B$ with $A, B \in\C$, $\xi\in\B$ and $\xi^2=-1$. Squaring and equating to $q$ gives:
\[
q^2 = A^2 - B^2 + 2 A B \xi = A + \xi B = q
\]
Equating coefficients of $1$ and $\xi$ (that is, the real and imaginary parts with respect to $\xi$) gives:
\begin{align*}
A &= A^2 - B^2\\
B &= 2 A B
\end{align*}
The second of these equations requires that $A = \frac{1}{2}$.
Making this substitution into the first equation gives:
$\frac{1}{2} = \frac{1}{4} - B^2$, from which $B^2 = -\frac{1}{4}$. Since $B$ is complex, the
only solutions are $B = \pm\I/2$. Thus $q = \frac{1}{2}\pm\frac{1}{2}\xi\I$.
\end{proof}
The roots of $-1$ in \B include two trivial cases, as noted in \cite[Theorem 2.1]{10.1007/s00006-006-0005-8}.
These are $\pm\I$ which yield the trivial idempotents $q=0$ and $q=1$;
and $\pm\mu$ where $\mu$ is a real pure quaternion root of $-1$,
which gives a subset of the idempotents of the form $q = \frac{1}{2}\pm\frac{1}{2}\mu\I$.

A simple example of an idempotent is given by choosing $\mu=\i$, thus:
$\left(\frac{1}{2} + \frac{1}{2}\I\i\right)^2 = \frac{1}{4} + \frac{1}{4}\I^2\i^2 + 2×\frac{1}{4}\I\i
= \frac{1}{2} + \frac{1}{2}\I\i$.

Clearly, since any idempotent is a solution of $q^2=q$, it is also a solution of $q(q-1)=0$.
Therefore, except for the trivial idempotents $0$ and $1$,
every biquaternion idempotent is also a divisor of zero.
In fact, as the following theorem shows, the connection between the idempotents and
the divisors of zero is much closer, since all the non-pure divisors of zero are
related to one of the idempotents by a complex scale factor.

\begin{theorem}
\label{nonpureidempotents}
All non-pure biquaternion divisors of zero are of the form $p = \alpha q$ where
$\alpha\in\C\mbox{\ and\ }\alpha\ne0$
and $q$ is a biquaternion idempotent as defined in Theorem~\ref{idempotent}.
\end{theorem}
\begin{proof}
The proof is by construction, from an arbitrary divisor of zero.
Let $p = A + \xi B$ with $A,B\in\C,\ \xi\in\B$ and $\xi^2=-1$, be an arbitrary biquaternion divisor of zero,
so that $\norm{p} = 0$. Since $\norm{p} = p\qconjugate{p}$ it is required that $(A + \xi B)(A - \xi B) = 0$.
Multiplying out, we have:
\[
(A+\xi B)(A-\xi B) = A^2 - \xi A B + \xi A B - \xi^2B^2 = A^2 - \xi^2B^2 =0
\]
Since $\xi$ is a root of $-1$, this yields the condition $A^2 + B^2 = 0$.
$A$ is required to be non-zero because $p$ is non-pure, so we require $A^2 = - B^2$,
which can only be satisfied if $A=\pm\I B$ (or $B = \pm\I A$).
Writing\footnote{If we write instead $A=\pm\I B$, then $\alpha=\pm 2\I B=2A$,
and the same result is reached after slightly more algebraic manipulation.}
$B = \pm\I A$, we have $p = A \pm\xi\I A = A(1\pm\xi\I)$.
Dividing by $\alpha=2A$, we obtain $p/\alpha = q$, where $q=\frac{1}{2}\pm\frac{1}{2}\xi I$,
an idempotent as in Theorem~\ref{idempotent}.
\end{proof}
\begin{corollary}
\label{dividem}
Dividing an arbitrary non-pure biquaternion divisor of zero by twice its scalar part always yields an idempotent.
\end{corollary}
Corollary~\ref{dividem} holds even when the divisor of zero is itself an idempotent,
since by Theorem~\ref{idempotent},
the scalar part has a value of $\frac{1}{2}$, and dividing by twice this value has no effect.
\begin{corollary}
Squaring an arbitrary non-pure divisor of zero, $p=\alpha q$, as in Theorem~\ref{nonpureidempotents}
yields $p^2 = \alpha p$, where $\alpha$ is twice the scalar part of $p$.
\end{corollary}
\begin{proof}
If $p = \alpha q$ and $q$ is idempotent, then $p^2 = \alpha^2 q^2 = \alpha^2 q = \alpha p$.
\end{proof}

In Theorem~\ref{vanish} we reviewed Hamilton's result that a divisor of zero must have
perpendicular real and imaginary parts with equal norm.
In Theorem~\ref{idempotent} we used a different representation to find the idempotents,
while in Theorem~\ref{nonpureidempotents} we proved that a non-pure
divisor of zero differs from an idempotent only by a complex scale factor.
It is instructive to verify that every idempotent satisfies
the two conditions required for a biquaternion to be a divisor of zero, which we now do.
\begin{proposition}
Every biquaternion idempotent $q = \frac{1}{2}\pm\frac{1}{2}\xi\I$ as in Theorem~\ref{idempotents}
satisfies the conditions in Theorem~\ref{vanish} and is therefore a divisor of zero.
\end{proposition}
\begin{proof}
We must show that the real and imaginary parts of $q$ have equal norm,
and that their inner product vanishes.
To do this we make use of the properties of $\xi$, which is a square root of $-1$.
These were given in \cite{arXiv:math.RA/0506190,10.1007/s00006-006-0005-8}:
$\xi$ must be a \emph{pure} biquaternion;
writing $\v\alpha=\Re(\xi)$ and $\v\beta=\Im(\xi)$,
we require that $\v\alpha\perp\v\beta$ and
$\norm{\v\alpha}-\norm{\v\beta}=1$.
Expressing $q$ in terms of the components of $\xi$ we have
$q = \frac{1}{2}\pm\frac{1}{2}(\v\alpha+\I\v\beta)\I = \frac{1}{2}\pm\frac{1}{2}(-\v\beta+\I\v\alpha)$.
Separating $q$ into real and imaginary parts, we have $\Re(q)=\frac{1}{2}(1\pm\v\beta)$
and $\Im(q)=\frac{1}{2}\v\alpha$.
The norms of the real and imaginary parts are as follows:
$\norm{\frac{1}{2}(1\pm\v\beta)}= \frac{1}{4}\norm{1\pm\v\beta} = \frac{1}{4}(1+\norm{\v\beta})$ and
$\norm{\frac{1}{2}\v\alpha}=\frac{1}{4}\norm{\v\alpha}$.
From the relation between the norms of $\v\alpha$ and $\v\beta$,
we have that $\norm{\v\alpha}=1+\norm{\v\beta}$,
and making this substitution,
we verify that the real and imaginary parts of $q$ do indeed have equal norms.
The inner product of the real and imaginary parts of $q$ is given by:
$\inner{\frac{1}{2}(1\pm\v\beta)}{\frac{1}{2}\v\alpha}=
\inner{\frac{1}{2}}{\frac{1}{2}\v\alpha}+
\inner{\pm\frac{1}{2}\v\beta}{\frac{1}{2}\v\alpha}$
which is zero because the inner product of the real scalar
$\frac{1}{2}$ with the vector $\frac{1}{2}\v\alpha$ is zero,
and $\v\alpha\perp\v\beta$.
\end{proof}

\section{Nilpotents}
\label{nilpotents}
The existence of biquaternion nilpotents was known to Hamilton,
and he realised that all such nilpotents were divisors of zero
\cite[Lecture VII, §\,674, pp\,671--3]{Hamilton:1853}. We
start by reviewing his result that a pure biquaternion divisor of zero has a vanishing square.
We then establish that all nilpotent biquaternions must be pure, and we show that the three complex components of the biquaternion have a null sum of squares.

\begin{theorem}
\label{vanishsquare}
\textnormal{\cite[Lecture VII, §\,672, p\,669]{Hamilton:1853}}\\
Let $p = p_r + \I p_i$ be a non-zero pure biquaternion with vanishing semi-norm.
Then the square of $p$ will be zero.
\end{theorem}
\begin{proof}
From Theorem \ref{vanish}, $\norm{p_r} = \norm{p_i}$, therefore we may divide by their common
norm, and obtain $p/\norm{p_r} = \mu + \I\nu$ where $\mu$ and $\nu$ are unit
pure quaternions which must be perpendicular because $\inner{p_r}{p_i}=0$.
Then:
\[
\left(\frac{p}{\norm{p_r}}\right)^2 = (\mu + \I\nu)^2 = \mu^2 - \nu^2 + \I(\mu\nu + \nu\mu)
\]
which vanishes because the squares of the two unit pure quaternions are $-1$, and the term on the right
is zero because the product of two perpendicular unit pure quaternions changes sign when the order of
the product is reversed.
\end{proof}
We can easily construct an example of a pure biquaternion with vanishing square to
demonstrate this result. Let $p = \i + \I\j$. Then we have
$p^2 = (\i + \I\j)(\i + \I\j) = \i^2 - \j^2 + \I(\i\j+\j\i) = 0 + \I(\k - \k) = 0$.

\begin{lemma}
\label{nilpotentsarepure}
All nilpotent biquaternions are pure.
\end{lemma}
\begin{proof}
Write an arbitrary biquaternion in the following form: $q = A + \xi B$, where $A, B\in\C$
and $\xi$ is a biquaternion root of $-1$ as in \cite{arXiv:math.RA/0506190,10.1007/s00006-006-0005-8}.
Then the square is:
\[
q^2 = A^2 + 2AB\xi + (\xi B)^2 = 0
\]
Now, $\xi B$ may be a nilpotent, in which case its square vanishes. Considering this
case first, we would have: $A^2 + 2AB\xi = 0$, and equating scalar and vector parts
of the left- and right-hand sides, we see directly that $A^2=0$ and hence $A=0$.

Now considering the alternative (that $\xi B$ is not a nilpotent),
since $\xi$ is a root of $-1$, $(\xi B)^2 = -B$, and we would have:
\[
A^2 - B^2 + 2AB\xi = 0
\]
Equating scalar and vector parts of the left- and right-hand sides, we obtain:
\begin{align*}
A^2 - B^2 &= 0\\
2AB\xi    &=0
\end{align*}
Since $\xi$ is not zero (because it is a root of $-1$),
the second equation can only be satisfied by $A=0$ and/or $B=0$.
Either choice would require the other to be zero, from the first equation.
Therefore, we conclude that $\xi B$ \emph{is} a nilpotent, and, since $A$ is zero,
$q$ is pure.
Since we imposed no restrictions on $q$ other
than those stated above, all nilpotent biquaternions must be pure.
\end{proof}
\begin{theorem}
\label{xyzsquares}
A biquaternion $q = W + X\i + Y\j + Z\k$ is nilpotent if $W = 0$ and $X^2 + Y^2 + Z^2 = 0$.
\end{theorem}
\begin{proof}
From Lemma~\ref{nilpotentsarepure}, only a pure biquaternion can be nilpotent and hence
$W$ must be zero.
Squaring $q$ gives the following result:
\[
q^2 =
\begin{array}{r@{}lr@{}l@{}lr@{}l@{}l}
- & X^2   & + & X & Y\k  & - & X & Z\j\\
- & X Y\k & - &   & Y^2  & + & Y & Z\i\\
+ & X Z\j & - &   & YZ\i & - &   & Z^2
\end{array}
\]
from which we can see that the `off-diagonal' terms cancel, leaving $q^2 = -(X^2 + Y^2 + Z^2)$.
\end{proof}
We have established in Lemma~\ref{nilpotentsarepure} that any nilpotent biquaternion must be pure.
In general, a pure biquaternion can be represented in
the form $\xi B$ where $B$ is complex and $\xi$ is a biquaternion root of $-1$,
as discussed in §\,\ref{notations}. In Theorem~\ref{xyzsquares} we have established that the sum of the squares of the complex
components of $\xi B$ must be zero.
This means that $B$ itself must be zero, since it is the square root of the semi-norm
($B=\sqrt{X^2+Y^2+Z^2}$), and the semi-norm of a divisor of zero vanishes.
Therefore, it is not possible to compute the axis, $\xi$, of a nilpotent biquaternion.

\begin{proposition}
\label{nilpotentsaredivisors}
Every biquaternion nilpotent $q = X\i + Y\j + Z\k$ with $X^2 + Y^2 + Z^2 = 0$
as in Theorem~\ref{xyzsquares} satisfies the conditions in Theorem~\ref{vanish}
and is therefore a divisor of zero.
\end{proposition}
\begin{proof}
The proof is as given in Proposition~\ref{normequivalence} with the additional
restriction $W=0$.
\end{proof}

\begin{corollary}
Every biquaternion nilpotent $p$ may be `normalised' to the form $\mu+\I\nu$
as in Theorem~\ref{vanishsquare}.
\end{corollary}
Even with the normalisation offered by Theorem~\ref{vanishsquare},
it remains impossible to find the axis,
since $\mu+\I\nu$ is still a nilpotent and therefore has a vanishing norm.
However, this step reduces the set of nilpotents to a well-understood set
and makes clear that the set of nilpotents has a simple structure based on
perpendicular unit pure quaternions.

\section{Composition of the biquaternion set of divisors of zero}
In Theorem~\ref{nonpureidempotents} we have shown that all non-pure divisors of zero
yield an idempotent on dividing by twice their scalar part (Corollary \ref{dividem});
and in Theorem \ref{vanishsquare} that all pure divisors of zero are nilpotent.
Therefore the set of zero divisors is composed of the set of idempotents multiplied
by a scalar (complex number) plus the set of nilpotents.


\end{document}